\documentclass[12pt,twoside]{article}
\usepackage{amsmath,amssymb,amsfonts}
\usepackage{amsthm}
\newtheorem{thm}{Theorem}[section]
\newtheorem{lem}[thm]{Lemma}

\newenvironment{pf}[1][Proof]{\noindent\textbf{#1.} }{\hfill \rule{0.5em}{0.5em}}

\DeclareMathOperator{\divg}{div}

\newcommand{\Real}{\mathbb{R}}

\newcommand{\abs}[1]{\left\vert#1\right\vert}
\newcommand{\norm}[1]{\left\Vert#1\right\Vert}

\title{Multiplicity of solutions for a class of quasilinear equations involving critical Orlicz-Sobolev nonlinear term}

\author{Jefferson A. Santos\footnote{J. A. Santos, was partially supported by PROCAD/CAPES-Brazil, e-mail: jefferson@dme.ufcg.edu.
br}\\
Universidade Federal de Campina Grande\\
Unidade Acad\^emica de Matem\'atica e Estat\'istica\\
CEP:58429-900, Campina Grande - PB, Brazil.}

\date{}

\begin{document}
\maketitle



\begin{abstract}
In this work, we study the existence and multiplicity of solutions
for a class of problems involving the $\phi$-Laplacian operator in a
bounded domain, where the nonlinearity has a critical growth. The
main tool used is the variational method combined with the genus
theory for even functionals.
\end{abstract}

\section{Introduction}
In this paper, we consider the existence and multiplicity of solutions for the following class of
quasilinear problem
\begin{align}
\left\{
\begin{array}{rcl}
  -\divg\left( \phi(\mid\nabla u\mid)\nabla u \right) & = & \lambda \phi_*(|u|)u + f(x, u), \quad \Omega\\
   u & =  & 0, \quad \partial \Omega  \\
\end{array}
\right.\tag{$P_{\lambda}$}\label{prob1}
\end{align}
where $ \Omega \subset \Real^{N} $ is a bounded domain in $\mathbb{R}^{N}$ with smooth boundary, $ \lambda $ is a positive parameter and $\phi : (0,+\infty) \to  \mathbb{R}$ is a continuous function verifying
$$
(\phi(t)t)'>0 \,\,\, \forall t>0.
\eqno{(\phi_{1})}
$$

\noindent There exist $l, m \in (1,N)$ such that
$$
l\leq \displaystyle\frac{\phi(|t|)t^{2}}{\Phi(t)}\leq m \,\,\,
\forall t \not= 0, \eqno{(\phi_{2})}
$$
where $\Phi(t)=\displaystyle\int^{|t|}_{0}\phi(s)s
\, ds$, $l\leq m < l^{*}$, $l^{*}=\displaystyle\frac{lN}{N-l}$ and
$m^{*}=\displaystyle\frac{mN}{N-m}$. Moreover, $\phi_*(t)t$ is such that Sobolev conjugate function $\Phi_*$ of $\Phi$ is its primitive, that is, $\Phi_*(t)=\int_0^{|t|}\phi_*(s)sds$.

Related to function $f:\overline{\Omega} \times \mathbb{R} \to
\mathbb{R}$, we assume that:
\begin{enumerate}
\item[$(f_1)$] $f \in C(\overline{\Omega} \times \mathbb{R},\mathbb{R})$ ,  odd with respect  $ t $ and
\begin{align*}
f(x, t) = o\left( \phi(\mid t\mid)\mid t\mid\right), \quad |t| \to 0 \mbox{ uniformly in } x;\\
f(x, t) = o\left( \phi_*(\mid t\mid)\mid t\mid \right), \quad |t|
\to +\infty \mbox{ uniformly in } x;
\end{align*}
\item[$(f_2)$] There is $\theta\in(m,l^*)$ such that $ F(x, t) \leq \frac{1}{\theta} f(x, t)t$, \,\, for all $ t > 0 $ and a.e. in $ \Omega $, where
$ F(x, t) = \int_{0}^{t}f(x, s)ds $.
\end{enumerate}

 There often arise the problem $(P_\lambda)$ associated by a nonhomogenous nonlinearities $\Phi$ in the fields of physics (see \cite{Fukagai3}), e.g.

$$
\begin{array}{l}
i) \mbox{ nonlinear elasticity:} \,\, \Phi(t)=(1+|t|^2)^{\gamma}-1 \,\,\, \mbox{for} \,\,\, \gamma \in (1, \frac{N}{N-2} ). \\
\mbox{}\\
ii) \mbox{ plasticity:} \,\,  \Phi(t)=|t|^{p}ln(1+|t|) \,\,\, \mbox{for} \,\,\, 1<p_0<p<N-1 \,\,\, \mbox{with}\\
 \displaystyle p_0=\frac{-1+\sqrt{1+4N}}{2}.\\
iii)\mbox{ generalized Newtonian fluids: } \Phi(t)=\int_0^ts^{1-\alpha}\left(\sinh^{-1}s\right)^\beta ds, \ 0\leq\alpha\leq 1,\\ \ \beta>0.
\end{array}
$$

Our main results is the following

\begin{thm} \label{T1}
Assume that $(\phi_1)-(\phi_2)$ and $(f_1)-(f_2)$ are satisfied. Then, there
exist a sequence $ \{ \lambda_{k} \} \subset (0,+\infty)$ with $
\lambda_{k}< \lambda_{k+1} $, such that, for $ \lambda \in
(\lambda_{k}, \lambda_{k+1}) $, problem~(\ref{prob1})  has at least
$ k $ pairs of nontrivial solutions.
\end{thm}

The main difficulty to prove Theorem \ref{T1} is related to the fact that the nonlinearity $f$ has a critical growth, because in this case, it is not clear that functional energy associated with problem $(P_\lambda)$ satisfies the well known $(PS)$ condition, once that the embedding $W^{1,\Phi}(\Omega) \hookrightarrow L_{\Phi_*}(\Omega)$ is not compact. To overcome this difficulty, we use a version of  the concentration compactness lemma due to Lions for  Orlicz-Sobolev space found in Fukagai, Ito and Narukawa  \cite{Fukagai1}. We would like to mention that  Theorem \ref{T1} improves the main result found
in \cite{xinmin}.

 We cite the papers of Alves and Barreiro \cite{Barreiro}, Alves, Gon\c calves and Santos \cite{Abrantes}, Bonano, Bisci and Radulescu \cite{BBR}, Cerny \cite{Cerny}, Cl\'ement, Garcia-Huidobro and Man\'asevich \cite{VGMS}, Donaldson \cite{Donaldson}, Fuchs and Li \cite{Fuchs1},  Fuchs and Osmolovski \cite{Fuchs2}, Fukagai, Ito and Narukawa \cite{Fukagai1,Fukagai2}, Gossez \cite{Gossez}, Mihailescu and Raduslescu \cite{MR1, MR2}, Mihailescu and Repovs \cite{MD}, Pohozaev \cite{Pohozaev}  and references therein, where quasilinear problems like $(P_\lambda)$ have been considered in bounded and unbounded domains of $\mathbb{R}^{N}$. In some those papers, the authors have mentioned that this class of problem arises in a lot of applications, such as, nonlinear elasticity, plasticity and non-Newtonian fluids.

This paper is organized in the following way: In Section~\ref{Orlicz-Sobolev}, we collect some preliminaries on Orlicz-Sobolev spaces that will be used throughout the paper, which can be found in \cite{adams}, \cite{adams2}, \cite{Donaldson2} and \cite{Oneill}. In Section 3, we recall an abstract theorem involving genus theory that will use in the proof of Theorem \ref{T1} and prove some technical lemmas, and Section 4 we prove Theorem \ref{T1}.

\vspace{0.5 cm}

\section{Preliminaries on Orlicz-Sobolev space }\label{Orlicz-Sobolev}
First of all, we recall that a continuous function $A:\mathbb{R} \to [0,+\infty)$ is a $N$-function if: \\

\noindent $(A1)$ \, $A$  is convex. \\

\noindent $(A2)$ \, $A(t)=0 \Leftrightarrow t=0$. \\

\noindent $(A3)$ \, $ \displaystyle \frac{A(t)}{t} \stackrel{t \to 0}{\longrightarrow} 0 \,\,\, \mbox{and} \,\,\, \frac{A(t)}{t} \stackrel{t \to +\infty }{\longrightarrow}+\infty $. \\

\noindent $(A4)$ \, $A$ is even.

\vspace{0.5 cm}
In what follows, we say that a $N$-function $A$ verifies the $\Delta_{2}$-condition if, there exists $t_0\geq0$ and $k>0$ such that
$$
A(2t)\leq kA(t) \,\,\,  t\geq t_0.
$$
This condition can be rewritten of the following way: For each $s >0$, there exists $M_s>0$ and $t_0\geq0$ such that
$$
A(st) \leq M_s A(t) \,\,\,t \geq t_0.    \eqno{(\Delta_2)}
$$

Fixed an open set $\Omega \subset \mathbb{R}^{N}$ and a N-function $A$ satisfying $\Delta_{2}$-condition, the space $L_{A}(\Omega)$ is the vectorial space  of the measurable functions $u: \Omega \to \mathbb{R}$ such that
$$
\displaystyle\int_{\Omega}A(u) < \infty.
$$
The space $L_{A}(\Omega)$ endowed with Luxemburg norm, that is, with the norm given by
$$
|u|_{A}= \inf \biggl\{\alpha >0: \int_{\Omega}A\Big(\frac{u}{\alpha}\Big)\leq 1\biggl\},
$$
is a Banach space. The complement function of $A$, denoted by $\widetilde{A}(s)$, is given by the Legendre transformation, that is
$$
\widetilde{A}(s)=\displaystyle\max_{t \geq 0}\{st -A(t)\} \ \ \mbox{for} \ \ s \geq 0.
$$
The functions $A$ and $\widetilde{A}$ are complementary each other. Moreover, we have the Young's inequality given by
\begin{equation} \label{D2}
st \leq A(t) + \widetilde{A}(s) \,\,\,\,\,\, \forall t,s \geq 0.
\end{equation}
Using the above inequality, it is possible to prove a H\"{o}lder type inequality, that is,
\begin{equation} \label{D3}
\biggl|\displaystyle\int_{\Omega}u v\biggl|\leq 2|u|_{A}|v|_{\widetilde{A}}, \,\,\, \forall \,\, u \in L_{A}(\Omega) \,\,\, \mbox{and} \,\,\,  v \in L_{\widetilde{A}}(\Omega).
\end{equation}
Another important function related to function $A$, it is the Sobolev's conjugate function $A_{*}$ of $A$ defined by
$$
A^{-1}_{*}(t)=\displaystyle\int^{t}_{0}\displaystyle\frac{A^{-1}(s)}{s^{(N+1)/N}} ds, \ \mbox{for } t>0.
$$

When $A(t)=|t|^{p}$ for $1< p<N $, we have $A_{*}(t)={p^{*}}^{p^{*}}|t|^{p^{*}}$, where $p^{*}=\frac{pN}{N-p}$.

Hereafter, we denote by $W^{1,A}_{0}(\Omega)$ the Orlicz-Sobolev space obtained by the completion of $C^{\infty}_{0}(\Omega)$ with respect to norm
$$
\|u\|=|\nabla u|_{A}+ |u|_{A}.
$$

An important property that we must detach is: If $A$ and $\widetilde{A} $ verifying $\Delta_2$-condition, the spaces $L_{A}(\Omega)$ and $W^{1,A}(\Omega)$ are reflexive and separable. Moreover, the $\Delta_2$-condition also implies that
\begin{equation} \label{CV0}
u_n \to u \,\,\, \mbox{in} \,\,\, L_{A}(\Omega) \Longleftrightarrow \int_{\Omega}A(|u_n-u|) \to 0
\end{equation}
and
\begin{equation} \label{CV1}
u_n \to u \,\,\, \mbox{in} \,\,\, W^{1,A}(\Omega) \Longleftrightarrow \int_{\Omega}A(|u_n-u|) \to 0 \,\,\, \mbox{and} \,\,\, \int_{\Omega}A(|\nabla u_n- \nabla u|) \to 0.
\end{equation}

Another important inequality  was proved by    Donaldson and Trudinger \cite{Donaldson}, which establishes that
for all open $\Omega \subset \mathbb{R}^{N}$ and there is a constant  $S_N=S(N) > 0$ such that
\begin{equation}\label{trudinger-emb}
\mid u\mid_{A_*}\leq S_N\mid\nabla
u\mid_{A}, \ u\in W_0^{1,A}(\Omega).
\end{equation}
Moreover, exist $ C_0>0$ such that
\begin{equation}\label{Poincare}
\int_\Omega A(u)\leq C_0\int_\Omega A(|\nabla u|), \ u\in W_0^{1,A}(\Omega).
\end{equation}
This inequality shows the below embedding is continuous
$$
W_0^{1,A}(\Omega) \stackrel{\hookrightarrow}{\mbox{\tiny cont}} L_{A_*}(\Omega).
$$
If bounded domain $\Omega$ and the limits below hold
\begin{equation} \label{M1}
\limsup_{t \to 0}\frac{B(t)}{A(t)}< +\infty \,\,\, \mbox{and} \,\,\, \limsup_{|t| \to +\infty}\frac{B(t)}{A_{*}(t)}=0,
\end{equation}
the embedding
\begin{equation} \label{M2}
W_0^{1,A}(\Omega) \hookrightarrow L_{B}(\Omega)
\end{equation}
is compact.

The next four lemmas involving the functions $\Phi, \widetilde{\Phi}$ and $\Phi_{*}$ and theirs proofs can be found in \cite{Fukagai1}. Hereafter, $\Phi$ is the $N$-function given in introduction and $\widetilde{\Phi},\Phi_{*}$ are the complement and conjugate functions of $\Phi$ respectively.
\begin{lem}
Assume $(\phi_1)-(\phi_2)$. Then,
$$
\Phi(t) = \int_0^{|t|} s \phi(s) ds, \,\,\,
$$
is a $N$-function with $\Phi, \widetilde{\Phi} \in \Delta_2$.  Hence, $L_\Phi(\Omega), W^{1,\Phi}(\Omega)$ and $W_0^{1,\Phi}(\Omega)$ are reflexive and separable spaces.
\end{lem}

\begin{lem} \label{F0} The functions $\Phi$, $\Phi_*$,  $\widetilde{\Phi}$ and $\widetilde{\Phi}_*$ satisfy the inequality
\begin{equation} \label{D1}
\widetilde{\Phi}(\phi(|t|)t) \leq \Phi(2t)  \mbox{ and } \widetilde{\Phi}_*(\phi_*(|t|)t) \leq \Phi_*(2t)\,\,\, \forall t \geq 0.
\end{equation}
\end{lem}

\begin{lem} \label{F1} Assume that $(\phi_1)-(\phi_2)$ hold and let $\xi_{0}(t)=\min\{t^{l},t^{m}\}$,\linebreak $ \xi_{1}(t)=\max\{t^{l},t^{m}\},$ for all $t\geq 0$. Then,
$$
\xi_{0}(\rho)\Phi(t) \leq \Phi(\rho t) \leq  \xi_{1}(\rho)\Phi(t) \;\;\; \mbox{for} \;\; \rho, t \geq 0
$$
and
$$
\xi_{0}(|u|_{\Phi}) \leq \int_{\Omega}\Phi(u) \leq \xi_{1}(|u|_{\Phi})  \;\;\; \mbox{for} \;\; u \in L_{\Phi}(\Omega).
$$
\end{lem}

\begin{lem} \label{F3} The function $\Phi_*$ satisfies the following inequality
$$
l^{*} \leq \frac{\Phi'_*(t)t}{\Phi_{*}(t)} \leq m^{*} \,\,\, \mbox{for} \,\,\, t > 0.
$$
\end{lem}
As an immediate consequence of the Lemma \ref{F3}, we have the following result

\begin{lem} \label{F2} Assume that $(\phi_1)-(\phi_2)$ hold and let  $\xi_{2}(t)=\min\{t^{l^{*}},t^{m^{*}}\},$ $\xi_{3}(t)=\max\{t^{l^{*}},t^{m^{*}}\}$ for all $t\geq 0$. Then,
$$
\xi_{2}(\rho)\Phi_*(t) \leq \Phi_*(\rho t) \leq \xi_{3}(\rho)\Phi_*(t) \;\;\; \mbox{for} \;\; \rho, t \geq 0
$$
and
$$
\xi_{2}(|u|_{\Phi_*}) \leq \int_{\Omega}\Phi_*(u)dx \leq \xi_{3}(|u|_{\Phi_*})  \;\;\; \mbox{for} \;\; u \in L_{A_*}(\Omega).
$$
\end{lem}
\begin{lem} \label{lem Phiest}
Let $\widetilde{\Phi}$ be the complement of $\Phi$ and put
$$
\xi_4(s)=\min\{s^{\frac{l}{l-1}}, s^{\frac{m}{m-1}}\}\ \mbox{and}\
\xi_5(s)=\max\{s^{\frac{l}{l-1}}, s^{\frac{m}{m-1}}\}, \ s\geq0.
$$
Then the following inequalities hold
$$
\xi_4(r)\widetilde{\Phi}(s)\leq\widetilde{\Phi}(r s)\leq
\xi_5(r)\widetilde{\Phi}(s),\ r,s\geq0
$$
and
$$
\xi_4(| u|_{\widetilde{\Phi}})\leq
\int_\Omega\widetilde{\Phi}(u)dx\leq \xi_5(|u|_{\widetilde{\Phi}}),\
u\in L_{\widetilde{\Phi}}(\Omega).
$$
\end{lem} \vskip0.5cm

\section{An abstract theorem and technical lemmas}

In  this section we recall an important abstract theorem involving genus theory, which will use in the proof of Theorem \ref{T1}. After, we prove some technical lemmas that will use to show that the energy functional associated with problem~(\ref{prob1}) verifies the hypotheses of the abstract theorem.

\subsection{An abstract theorem}
Let $ E $ be a real Banach space and $ \Sigma $  the family of sets $ Y \subset E \backslash \{ 0 \} $ such that $ Y $ is closed in $ E $ and symmetric with respect to 0, that is,
\[
\Sigma =  \left\{ Y \subset E \backslash \{ 0 \}; Y \mbox{ is closed in } E \mbox{ and } Y = -Y \right\}.
\]
Hereafter, let us denote by $ \gamma (Y) $ the genus of $ Y \in \Sigma$ \, ( see \cite[pp. 45]{Rabinowitz1} ). Moreover, we set
\[
K_{c} = \left\{ u \in E; I(u) = c \mbox{ and } I'(u) = 0 \right\}
\]
and
\[
A_{c} = \left\{ u \in E; I(u) \leq c  \right\}.
\]

Next, we recall a version of the Mountain Pass Theorem for even functional. For details of the proof,  see \cite{Rabinowitz1}.

\begin{thm}\label{genus}
Let $ E $ be an infinite dimensional Banach space with \linebreak $ E = V \oplus X $, where $ V $ is finite dimensional and let $ I \in C^{1}(E,\mathbb{R}) $ be a even function with $ I(0) = 0 $ and satisfying :
\begin{enumerate}
\item [($I_{1}$)] there are constants $ \beta, \rho > 0 $ such that $ I(u) \geq \beta > 0 $, for each $ u \in \partial B_{\rho} \cap X $;
\item [($I_{2}$)] there is $ \Upsilon > 0 $ such that $ I $ satisfies the $(PS)_{c}$ condition, for $ 0 < c < \Upsilon$;
\item [($I_{3}$)] for each finite dimensional subespace $ \widetilde{E} \subset E $, there is  $ R  = R(\widetilde{E})>0 $ such that $ I(u) \leq 0 $ for all $ x \in  \widetilde{E} \backslash B_{R}(0)$.
\end{enumerate}

Suppose $ V $  is $ k$ dimensional and $ V = span\{e_{1}, \cdots , e_{k}\}$. For $ m \geq k $, inductively choose $ e_{m+1} \not\in E_{m} := span\{e_{1}, \cdots , e_{m}\}$. Let $ R_{m} = R(E_{m}) $ and $ D_{m} = B_{R_{m}} \cap E_{m} $. Define
\begin{eqnarray}
G_{m} := \left\{ h \in C(D_{m}, E); h \text{ is odd and } h(u) = u, \forall u \in \partial B_{R_{m}} \cap E_{m} \right\}
\end{eqnarray}
and
\begin{eqnarray}
\Gamma_{j} := \left\{ h(\overline{D_{m} \backslash Y}) ; h \in G_{m}, m \geq j, Y \in \Sigma, \text{ and } \gamma(Y) \leq m - j\right\}.
\end{eqnarray}
For each $ j \in \mathbb{N} $, let
\begin{eqnarray}
c_{j} = \inf_{K \in \Gamma_{j}} \max_{u \in K} I(u).
\end{eqnarray}
Then, $ 0 < \beta \leq c_{j} \leq c_{j+1}$ for $ j > k $, and if $ j > k,  c_{j} < \Upsilon$ and $ c_{j} $ is critical value of $ I $. Moreover, if $ c_{j}=c_{j+1}= \cdots = c_{j+l} = c < \Upsilon $ for $ j > k $, then $ \gamma(K_{c}) \geq l + 1 $.
\end{thm}

\subsection{Technical lemmas}

Associated with the problem~(\ref{prob1}), we have the energy functional \linebreak $ J_{\lambda} : W^{1,\Phi}_{0}(\Omega) \to \Real $ defined by
$$
J_{\lambda} (u) = \int_{\Omega} \Phi(\mid\nabla u\mid)-\lambda \int_{\Omega} \Phi_*(u)  - \int_{\Omega} F(x, u).
$$

By conditions $(f_1)-(f_2)$, $ J_{\lambda} \in C^{1} \left( W^{1,\Phi}_{0}(\Omega), \Real \right) $ with
\begin{align}
J'_{\lambda}(u)\cdot v = & \int_{\Omega} \phi(\mid\nabla u\mid) \nabla u \nabla v  - \lambda \int_{\Omega} \phi_*(\mid u\mid)u v  - \int_{\Omega} f(x, u)v, \nonumber
\end{align}
for any $ u, v \in W^{1,\Phi}_{0}(\Omega) $. Thus,  critical points of $ J_{\lambda} $ are weak solutions of  problem~(\ref{prob1}).

\begin{lem}
Under the conditions $(f_1) -(f_2) $, $ J_{\lambda} $ satisfies ($I_{1}$).
\end{lem}
\begin{pf}

On the other hand, from $(f_1)-(f_2)$, given $\epsilon >0$, there exists $C_\epsilon >0$ such that
\begin{align}
\abs{F(x, t)} \leq \epsilon\Phi(t)+C_{\epsilon}\Phi_*(t),\ \quad \forall (x, t) \in \bar{\Omega} \times \Real.\label{cond_G}
\end{align}

\noindent Combining \eqref{Poincare} with (\ref{cond_G}),
$$
J_{\lambda}(u) \geq \left(1-\epsilon C_0\right)\int_{\Omega}\Phi(\mid\nabla u\mid)  - (1+C_{\epsilon}) \int_{\Omega} \Phi_*(u).
$$

\noindent For $\epsilon$ is small enough and $ \norm{u} = \rho \simeq 0 $, it follows from (\ref{trudinger-emb}) and
Lemma \ref{F3}
\[
J_{\lambda}(u) \geq C_1\mid\nabla u\mid_\Phi^m-C_2S_N^{l^*}\mid \nabla u\mid_\Phi^{l^*}.
\]
for some positive constants $C_1$ and $C_2$. Once that, $ m < l^* $, if $\rho$ is small enough, there is $\beta >0$ such that
\[
J_{\lambda}(u) \geq \beta > 0  \,\,\, \forall u \in \partial B_{\rho}(0),
\]
finishing the proof.
\end{pf}

\vspace{0.5 cm}

\begin{lem}
Under the conditions $(f_1)-(f_2) $, $ J_{\lambda} $ satisfies $(I_{3})$.
\end{lem}
\begin{pf}
Suppose  $(I_{3})$ does not hold. Then, there is a finite dimensional subspace $ \widetilde{E} \subset W^{1,\Phi}_{0}(\Omega) $ and a sequence $ (u_{n}) \subset  \widetilde{E} \backslash B_{n}(0) $ verifying:
\begin{equation} \label{NOVAEQUA}
J_{\lambda}(u_{n}) > 0, \quad \forall n \in \mathbb{N} .
\end{equation}
A direct computation shows that given $\epsilon >0$, there is a constant $ M > 0 $ such that
\begin{align}\label{desM}
F(x,t) \geq -M - \epsilon \Phi_*(t)\,\,\, \forall (x,t) \in \bar{\Omega} \times \Real.
\end{align}
Consequently,
$$
J_{\lambda} (u_{n}) \leq  \int_{\Omega}\Phi(\abs{\nabla u_n})dx  - \lambda \int_{\Omega} \Phi_*(u_n)+ \epsilon\int_{\Omega}\Phi_*(u_n) + M \abs{\Omega}.
$$
Fixing $ \epsilon = \frac{\lambda}{2} $, and using Lemma \ref{F2}, we get
\begin{align}
J_{\lambda}(u_{n}) \leq \int_{\Omega}\Phi(\abs{\nabla u_{n}}) -  \frac{\lambda}{2}\xi_3(\abs{u_{n}}_{\Phi_*})  + M \abs{\Omega}.
\end{align}
Once that $ \dim{\tilde{E}} < \infty $, we know that any two norms in $ \tilde{E} $ are equivalent. Then, using that $ \norm{u_{n}} \to \infty  $, we can assume that                $ |u_{n}|_{\Phi_*}> 1 $. Thereby, from  Lemmas \ref{F1} and \ref{F2},
$$
J_{\lambda}(u_{n}) \leq |\nabla u_n|_{\Phi}^m -\frac{\lambda}{2}|u_n|^{l^*}_{\Phi_*} + M\abs{\Omega}.
$$
Using again the equivalence of the norms in $ \tilde{E} $, there is $C>0$ such that
$$
J_{\lambda}(u_{n}) \leq  \parallel u_n\parallel^m - \frac{\lambda}{2}C\parallel u_n\parallel^{l^*} + M\abs{\Omega}.
$$
Recalling that $ m < l^*$, the above inequality implies that there is $n_0 \in \mathbb{N}$ such that
\begin{align*}
J_{\lambda}(u_{n}) < 0, \,\,\, \forall  n \geq n_0,
\end{align*}
which contradicts \eqref{NOVAEQUA}.
\end{pf}

\vspace{0.5 cm}

\vspace{0.5 cm}

\begin{lem}\label{limitada}
Under the conditions $(f_1) - (f_2)$, any $(PS)$ sequence for $J_\lambda$ is bounded in $W_0^{1,\Phi}(\Omega)$.
\end{lem}

\begin{pf} \, Let $ \{ u_{n} \} $ be a $(PS)_{d} $ sequence of  $ J_{\lambda} $. Then,
\[
J_{\lambda}(u_{n}) \to d \mbox{ and } J'_{\lambda}(u_{n}) \to 0 \mbox{ as} \,\,\, n \to +\infty.
\]
We claim that $ \{ u_{n} \} $ is bounded. Indeed, note that
\begin{align*}
J_{\lambda}(u_{n}) - \frac{1}{\theta}J'_{\lambda}(u_{n})u_{n} = &
 \int_{\Omega}\Phi(\abs{\nabla u_{n}})-\frac{1}{\theta}\int_\Omega\phi(\mid\nabla u_n\mid)\mid\nabla u_n\mid^2   \\
& -\lambda\int_\Omega\Phi_*(u_n)+\frac{\lambda}{\theta}\int_\Omega \phi_*(|u_n|)u_n^2\\
&-\int_\Omega F(x,u_n)+\frac{1}{\theta}\int_\Omega f(x,u_n)u_n.
\end{align*}
Consequently,
\begin{align*}
\lambda \int_{\Omega}\left(\frac{1}{\theta}\phi_*(|u_n|)u_n^2-\Phi_*(u_n)\right)= & J_{\lambda}(u_{n}) - \frac{1}{\theta}J'_{\lambda}(u_{n})u_{n} - \int_{\Omega}\Phi(\abs{\nabla u_{n}})  \\
 & +\frac{1}{\theta}\int_\Omega\phi(\mid\nabla u_n\mid)\mid\nabla u_n\mid^2  \\
 &+\int_{\Omega}\left( F(x, u_{n}) - \frac{1}{\theta}f(x, u_{n})u_{n} \right).
\end{align*}
Then, by $(\phi_2)$, $(f_2)$ and Lemma \ref{F3}, for $ n $ sufficiently large
$$
\lambda \left(\frac{l}{\theta}-1\right)\int_{\Omega}\Phi_*(u_n)  \leq C+1+\parallel u_n\parallel+\left(\frac{m}{\theta}-1\right)\int_\Omega\Phi(\mid\nabla u_n\mid),
$$
which implies that
\begin{align*}
\left[\lambda\left(\frac{l}{\theta}-1\right)\right] \int_{\Omega}\Phi_*(u_n) \leq C +  \norm{u_{n}},
\end{align*}
where $ C$ is a positive constant, and so
\begin{align}\label{desestr}
\int_{\Omega}\Phi_*(u_n) dx \leq C \left( 1 + \norm{u_{n}} \right).
\end{align}
By (\ref{desM}) and (\ref{desestr})
\begin{align*}
 \int_{\Omega}\Phi(\mid\nabla u_n\mid)& \leq  J_{\lambda}(u_n) + \lambda \int_{\Omega}\Phi_*(u_n) +\int_{\Omega} G(x, u_{n}) dx \\
& \leq C + o_n(1) + (\lambda+\epsilon)\int_\Omega\Phi_*(u_n)\\
& \leq C(1+\parallel u_n\parallel)+o_n(1).\\
\end{align*}
Therefore, for $ n $ sufficiently large
$$
\int_{\Omega}\Phi(\mid\nabla u_n\mid) \leq C\left(1 + \norm{u_{n}} \right).
$$
If $ \norm{ u_{n}} > 1 $, it follows from Lemma \ref{F2}
\begin{align*}
\parallel u_{n}\parallel^{l} \leq C\left(1 + \norm{u_{n}} \right).
\end{align*}
Once that $ l > 1 $, the above inequality gives  that $ \{ u_{n}\} $ is bounded in $ W_{0}^{1, \Phi}(\Omega) $.

\end{pf}

As a consequence of the last result, if $\{u_n\}$ is a $(PS)$ sequence for $J_\lambda$, we can extract a subsequence of $\{u_n\}$, still denoted by $ \{u_{n}\} $ and $u \in W_0^{1,\Phi}(\Omega)$, such that
\begin{itemize}
  \item $ u_{n} \rightharpoonup u $ in $ W_{0}^{1, \Phi}(\Omega) $;
  \item $ u_{n} \rightharpoonup u $ in $ L_{\Phi_*}(\Omega) $;
  \item $ u_{n} \to u $ in $ L_{\Phi}(\Omega) $;
  \item $u_n(x) \to u(x)$ a.e. in $\Omega$;
\end{itemize}

From the concentration compactness lemma of Lions in Orlicz-Sobolev space found in \cite{Fukagai1}, there exist two nonnegative measures $ \mu, \nu \in \mathcal{M}(\Omega) $, a countable set $\mathcal{J}$, points $ \{ x_{j} \}_{j \in \mathcal{J}} $ in $ \Omega $ and sequences  $ \{ \mu_{j}\}_{j \in \mathcal{J}}, \{ \nu _{j}\}_{j \in \mathcal{J}} \subset [0, +\infty)$, such that
\begin{align}
\Phi(\abs{\nabla u_{n}}) \to \mu \geq \Phi(\abs{\nabla u}) + \sum_{j \in \mathcal{J}}\mu_{j}\delta_{x_{j}} \mbox{ in } \mathcal{M}( \Omega)\\
\Phi_*(u_{n}) \to \nu =  \Phi_*(u) + \sum_{j \in \mathcal{J}}\nu_{j}\delta_{x_{j}} \mbox{ in } \mathcal{M}( \Omega)
\end{align}
and
\begin{align}
\nu_j\leq\max\{S_N^{l^*}\mu_j^{\frac{l^*}{l}},S_N^{m^*}\mu_j^{\frac{m^*}{l}},S_N^{l^*}\mu_j^{\frac{l^*}{m}},
S_N^{m^*}\mu_j^{\frac{m^*}{m}}\},
\end{align}
where $S_N$ verifies (\ref{trudinger-emb}).

Next, we will show an important estimate from below for $\{\nu_i\}$. We have to prove a technical lemma.

\vspace{0.5 cm}

\begin{lem} \label{NU}
Under the conditions of Lemma~\ref{limitada}. If $\{u_n\}$ is a $(PS)$ sequence for $J_\lambda$ and $\{\nu_j\}$ as above,  then for each $j \in \mathcal{J}$,
$$
\nu_j \geq\left(\frac{l}{\lambda m^*}\right)^{\frac{\beta}{\beta-1}}S_N^{-\frac{\alpha}{\beta-1}} \,\,\, \mbox{or} \,\,\,  \nu_j = 0,
$$
for some $\alpha \in\{l^*,m^*\}$ and $\beta\in \{\frac{l^*}{l},\frac{m^*}{l},\frac{l^*}{m},\frac{m^*}{m}\}$.
\end{lem}
\begin{pf} Let $ \psi \in C^{\infty}_{0}(\mathbb{R}^{N}) $ such that
$$
 \psi (x)= 1  \,\,\, \mbox{in} \,\,\,  B_{\frac{1}{2}}(0) ,  \,\,\, supp\,\psi \subset B_{1}(0) \,\, \mbox{and} \,\,\,  0 \leq \psi(x) \leq 1 \,\,\, \forall x \in \mathbb{R}^{N}.
$$
 For each $j\in \Gamma$ and $ \epsilon > 0 $, let us define
$$
\psi_{\epsilon}(x) = \psi\left(\frac{x-x_j}{\epsilon} \right), \,\,\, \forall  x \in \Real^{N}.
$$
Then $ \left\{ \psi_{\epsilon} u_{n} \right\} $ is bounded in $ W_{0}^{1, \Phi}(\Omega) $.  Since $J'_\lambda(u_n)\to 0$, we have
$$
J'_\lambda(u_n)(\psi_{\epsilon}u_n)=o_n(1),
$$
or equivalently,
\begin{align}
  \int_{\Omega}\phi( |\nabla u_{n}|)\nabla u_n\nabla\left(u_n \psi_{\epsilon}\right) &=o_n(1)+  \lambda \int_{\Omega} \phi_*(|u_{n}|)u_n^2 \psi_{\epsilon} +      \int_{\Omega} f(x,u_{n})u_{n} \psi_{\epsilon}\label{der_nabla_un_phi}\nonumber\\
&\leq o_n(1)+\lambda m^*\int_\Omega \Phi_*(u_n)\psi_\epsilon+\int_\Omega f(x,u_n)u_n \psi_\epsilon.
\end{align}

By the compactness Lemma of Strauss \cite{chabro}
\begin{equation}\label{conv_grad3}
\lim_{n \to \infty} \int_{\Omega} f(x, u_{n})u_{n}\psi_{\epsilon} = \int_{\Omega} f(x, u)u\psi_{\epsilon}.
\end{equation}

On the other hand, by $(\phi_2)$
\begin{align}\label{eq est03}
  \int_{\Omega}\phi( |\nabla u_{n}|)\nabla u_n\nabla\left(u_n \psi_{\epsilon}\right)& =\int_{\Omega}\phi( |\nabla u_{n}|)|\nabla u_n|^2 \psi_{\epsilon} + \int_{\Omega}\phi( |\nabla u_{n}|)\left(\nabla u_n\nabla \psi_{\epsilon}\right)u_n\nonumber\\
&\ \ \geq l\int_\Omega \Phi(|\nabla u|) \psi_{\epsilon}+ \int_{\Omega}\phi( |\nabla u_{n}|)\left(\nabla u_n\nabla \psi_{\epsilon}\right)u_n.
\end{align}
By Lemmas \ref{F0} and \ref{lem Phiest}, the sequence $\{|\phi(|\nabla u_n|)\nabla u_n|_{\widetilde{\Phi}}\}$ is bounded. Thus, there is a subsequence $\{u_n\}$ such that
$$
\phi(|\nabla u_n|)\nabla u_n\rightharpoonup \widetilde{w}_1 \mbox{ weakly in }L_{\widetilde{\Phi}}(\Omega,\mathbb{R}^N),
$$
for some $\widetilde{w}_1\in L_{\widetilde{\Phi}}(\Omega,\mathbb{R}^N)$. Since  $u_n \to u$ in $L_\Phi( \Omega)$,
$$
\int_\Omega \phi(|\nabla u_n|)(\nabla u_n\nabla \psi_\epsilon)u_n\to \int_\Omega (\widetilde{w}_1\nabla\psi_\epsilon)u.
$$
Thus, combining (\ref{der_nabla_un_phi}), (\ref{conv_grad3}), (\ref{eq est03}) and letting $n\to \infty$, we have
\begin{equation}\label{eq est02}
l\int_\Omega \psi_\epsilon d\mu+\int_\Omega (\widetilde{w}_1\nabla
\psi_\epsilon)u \leq \lambda m^*\int_\Omega \psi_\epsilon d\nu
+\int_\Omega f(x,u)u\psi_\epsilon.
\end{equation}
Now we show that the second term of the left-hand side converges 0
as $\epsilon \to 0.$
First, we show the claim:\\
{\bf Claim 1:} $\{f(x,u_n)\}$ is bounded in $L_{\widetilde{\Phi}_*}(\Omega)$.\\
In fact, by $(f_1)$ and Lemma \ref{F0} we have
\begin{align*}
\int_\Omega\widetilde{\Phi}_*( f(x,u_n))&\leq c_1\int_\Omega \widetilde{\Phi}_*(\phi_*(|u_n|)u_n)
+c_2\int_{[|u_n|>1]}\widetilde{\Phi}_*(\phi(|u_n|)u_n)\\
&+c_3\int_{[|u_n|\leq 1]}\widetilde{\Phi}_*(\phi(|u_n|)u_n)\\
&\leq c_1\int_\Omega\Phi_*(u_n)+c_2\int_{[|u_n|>1]}\widetilde{\Phi}_*(\phi(|u_n|)u_n)+c_3|\Omega|.
\end{align*}
Hence, by $(\phi_2)$, Lemma \ref{F1} and $m<l^*$,
\begin{align*}
\int_\Omega\widetilde{\Phi}_*( f(x,u_n))&\leq C_1\int_\Omega \widetilde{\Phi}_*(\phi_*(|u_n|)u_n)+C_2\int_{[|u_n|>1]}\widetilde{\Phi}_*(|u_n|^{m-1})+C_3|\Omega|\\
&\leq C_1\int_\Omega\Phi_*(u_n)+C_2\int_{[|u_n|>1]}\widetilde{\Phi}_*(|u_n|^{l^*-1})+C_3|\Omega|.
\end{align*}
Now, by Lemmas \ref{F3} and \ref{F0}
$$
\int_\Omega\widetilde{\Phi}_*( f(x,u_n))\leq K_1\int_\Omega \Phi_*(u_n)+K_2|\Omega|<+\infty.
$$
From Claim 1, there is a subsequence $\{u_n\}$ such that
$$
\phi_*(|u_n|)u_n+f(x,u_n)\rightharpoonup \widetilde{w}_2\mbox{ weakly in }L_{\widetilde{\Phi}_*}(\Omega),
$$
for some $\widetilde{w}_2\in L_{\widetilde{\Phi}_*}(\Omega)$. Since
\begin{align*}
J'_\lambda(u_n)v=&\int_\Omega\phi(|\nabla u_n|)\nabla u_n\nabla v-\int_\Omega\left(\phi_*(|u_n|)u_n+f(x,u_n)\right)v\\
&\to 0,
\end{align*}
as $n\to \infty$ for any $v\in W_0^{1,\Phi}(\Omega)$,
$$
\int_\Omega(\widetilde{w}_1\nabla v-\widetilde{w}_2v)=0,
$$
for any $v\in W_0^{1,\Phi}(\Omega)$. Substituting $v=u\psi_\epsilon$ we have
$$
\int_\Omega(\widetilde{w}_1\nabla(u\psi_\epsilon)-\widetilde{w}_2u\psi_\epsilon)=0.
$$
Namely,
$$
\int_\Omega(\widetilde{w}_1\nabla\psi_\epsilon)u=-\int_\Omega(\widetilde{w}_1\nabla u-\widetilde{w}_2u)\psi_\epsilon.
$$
Noting $\widetilde{w}_1\nabla u-\widetilde{w}_2u\in L^1(\Omega)$, we see that right-hand side tends to 0 as $\epsilon\to 0$. Hence we have
$$
\int_\Omega (\widetilde{w}_1\nabla \psi_\epsilon)u\to 0,
$$
as $\epsilon\to 0$.

Letting $\epsilon\to 0$ in (\ref{eq est02}), we obtain
\begin{align*}
l\mu_{j} \leq \lambda m^*\nu_{j}.
\end{align*}
Hence,
\[
S_N^{-\alpha}\nu_j \leq \mu_{j}^\beta\leq  \left(\frac{l\lambda}{m^*}\right)^{\beta} \nu_{j}^\beta,
\]
for some $\alpha\in\{l^*,m^*\}$, $\beta\in\{\frac{l^*}{l},\frac{m^*}{l},\frac{l^*}{m},\frac{m^*}{m}\}$, and so
$$
\nu_{j} \geq \left(\frac{l}{\lambda m^*}\right)^{\frac{\beta}{\beta-1}}S_{N}^{-\frac{\alpha}{\beta-1}} \quad \mbox{ or} \quad  \nu_{j} = 0 .
$$

\end{pf}

\begin{lem}
Assume that $(f_1)-(f_2) $. Then, $ J_{\lambda} $ satisfies $(PS)_{d} $  for $d\in(0,d_\lambda)$ where
$$
d_\lambda= \min\left\{\frac{l^*-\theta}{\theta S_N^{\frac{\alpha}{\beta-1}}\lambda^{\frac{1}{\beta-1}}}\left(\frac{l}{m^*}\right)^{\frac{\beta}{\beta-1}};\alpha\in\{l^*,m^*\},\ \beta\in\left\{\frac{l^*}{l},\frac{m^*}{l},\frac{l^*}{m},\frac{m^*}{m}\right\}\right\}.
$$
\end{lem}
\begin{pf} Once that
$$
 J_\lambda(u_{n}) = d + o_{n}(1) \,\,\, \mbox{and} \,\,\, J_\lambda'(u_{n})=o_n(1),
$$
 \begin{align}
d = \lim_{n \to \infty} I(u_{n}) &=   \lim_{n \to \infty} \left( J_\lambda(u_{n}) - \frac{1}{\theta} J_\lambda'(u_{n})u_{n}  \right) \nonumber\\
    &\geq \lim_{n \to \infty} \left[  \left( 1 - \frac{m}{\theta} \right)\int_{\Omega}\Phi( |\nabla u_{n}|)+ \lambda\left(\frac{l^*}{\theta}-1\right)\int_{\Omega} \Phi_*( u_{n}) \right. \nonumber\\
     &\left. - \int_{\Omega}\left( F(x, u_{n}) - \frac{1}{\theta} f(x,u_{n})u_{n} \right) \right]\nonumber\\
&\geq\lambda\left(\frac{l^*}{\theta}-1\right)\int_{\Omega} \Phi_*( u_{n}).
\end{align}
Recalling that
$$
\lim_{n \to \infty} \int_{\Omega} \Phi_*(u_n)dx =  \left[ \int_{\Omega} \Phi_*(u) + \sum_{j \in \mathcal{J}}\nu_{j} \right] \geq  \nu_{j},
$$
we derive that
\begin{align}
d \geq& \lambda \left(\frac{l^*}{\theta}-1\right)\left(\frac{l}{\lambda m^*}\right)^\frac{\beta}{\beta-1}S_N^{-\frac{\alpha}{\beta-1}}\nonumber\\
&=\left(\frac{l^*-\theta}{\theta}\right)\left(\frac{l}{m^*}\right)^{\frac{\beta}{\beta-1}}S_N^{\frac{\alpha}{\beta-1}}\lambda^{\frac{1}{1-\beta}},
\end{align}
for some $\alpha\in\{l^*,m^*\}$, $\beta\in\{\frac{l^*}{l},\frac{m^*}{l},\frac{l^*}{m},\frac{m^*}{m}\}$, which is an absurd. From this, we must have $ \nu_{j} = 0 $ for any $ j \in \mathcal{J}$, leading to
\begin{align}
\int_{\Omega}\Phi_*(u_{n}) \to \int_{\Omega}\Phi_*(u) \label{conv_critica}.
\end{align}
Combining the last limit with Br\'ezis and Lieb \cite{Brezis_Lieb}, we obtain
\[
\int_{\Omega} \Phi_*(u_{n} - u) \to 0 \mbox{ as} \quad n \to \infty,
\]
from where it follows by Lemma \ref{F2}
\begin{align}
u_{n} \to u \mbox{ in } L_{\Phi_*}(\Omega).
\end{align}
Now, as $ J_\lambda'(u_{n})u_{n} = o_{n}(1) $, the last limit gives
\begin{align*}
\int_{\Omega}\phi(|\nabla u_n|)\mid u_n\mid^2  =& \lambda \int_{\Omega}\phi_*( |u_n|)u_n^2  + \int_{\Omega} f(x, u_n)u_{n}+ o_n(1).
\end{align*}
In what follows, let us denote by $\{P_n\}$ the following sequence,
$$
P_{n}(x) = \langle \phi(|\nabla u_{n}(x)|) \nabla u_{n}(x) - \phi(|\nabla u(x)|) \nabla u(x), \nabla u_{n}(x) - \nabla u (x) \rangle.
$$
Since $\Phi$ is convex in $\mathbb{R}$ and $\Phi(|.|)$ is $C^1$ class in $\mathbb{R}^N$,  has $P_n(x) \geq 0. $
From definition of $\{P_n\}$,
\[
\int_{\Omega} P_{n} = \int_{\Omega} \phi(|\nabla u_{n}|)\mid\nabla u_n\mid^2 - \int_{\Omega}\phi( |\nabla u_{n}|) \nabla u_{n} \nabla u - \int_{\Omega} \phi(|\nabla u|) \nabla u \nabla(u_{n} - u).
\]
Recalling that $u_n \rightharpoonup u$ in $W_0^{1,\Phi}(\Omega)$, we have
\begin{align}
\int_{\Omega} \phi(|\nabla u|) \nabla u \nabla(u_{n} - u)  \to 0 \quad \mbox{ as } n \to \infty,
\end{align}
which implies that
\[
\int_{\Omega} P_{n} = \int_{\Omega} \phi(|\nabla u_{n}|) |\nabla u_{n}|^2- \int_{\Omega} \phi(|\nabla u_{n}|)\nabla u_{n} \nabla u + o_n(1).
\]
On the other hand, from $ J_\lambda'(u_{n})u_{n} = o_n(1)$ and $ J_\lambda'(u_{n})u = o_n(1)$, we derive
\begin{align*}
0\leq\int_{\Omega} P_{n} =& \lambda \int_{\Omega} \phi_*(|u_{n}|)|u_n|^2  - \lambda \int_{\Omega}\phi_*(|u_{n}|)u_{n} u\\
    &  +\int_{\Omega} f(x, u_{n})u_{n} - \int_{\Omega}f(x, u_{n}) u + o_n(1).
\end{align*}
Combining (\ref{conv_critica}) with the compactness Lemma of Strauss \cite{chabro}, we deduce that
\[
\int_{\Omega} P_{n} \to 0 \quad \mbox{as } n \to \infty.
\]
Applying a result due to Dal Maso and Murat \cite{Maso}, have that
\begin{equation} \label{E2}
u_{n} \to u \mbox{ in } W^{1,\Phi}_{0}(\Omega) .
\end{equation}
\end{pf}

\vspace{0.5 cm}

The next lemma is similar to \cite[Lemma~5]{xinmin} and its proof will be omitted.
\begin{lem}
Under the conditions $(f_1)$--$(f_2)$, there is sequence \linebreak $ \{M_{m}\} \subset (0,+\infty) $ independent of  $ \lambda $ with $ M_{m} \leq M_{m+1}$,  such that for any $ \lambda > 0 $
\begin{eqnarray}
c_{m}^{\lambda} = \inf_{K \in \Gamma_{m}} \max_{u \in K} J_{\lambda}(u) < M_{m}.
\end{eqnarray}

\end{lem}

\vspace{0.2 cm}

\section{Proof of Theorem \ref{T1}}

For each $ k \in \mathbb{N} $, choose $ \lambda_{k+1} $ such that
$$
M_{k} < d_{\lambda_k}.
$$
Thus, for $ \lambda \in (\lambda_{k}, \lambda_{k+1}) $,
\[
0 < c^{\lambda}_{1} \leq c^{\lambda}_{2} \leq \cdots  \leq  {c}^{\lambda}_{k} < M_{k} \leq d_\lambda.
\]
By Theorem~\ref{genus}, the levels $ c^{\lambda}_{1} \leq c^{\lambda}_{2} \leq \cdots  \leq  c^{\lambda}_{k} $ are critical values of $ J_{\lambda}$. Thus, if
$$
c^{\lambda}_{1} < c^{\lambda}_{2} < \cdots  <  c^{\lambda}_{k} ,
$$
functional $J_{\lambda}$ has at least $k$ critical points. Now, if $
c^{\lambda}_{j} = c^{\lambda}_{j+ 1} $ for some $ j = 1, 2, \cdots ,
k $, it follows from Theorem~\ref{genus} that $K_{c^{\lambda}_{j}} $
is an infinite set \cite[Cap. 7]{Rabinowitz1}. Then, in this case,
problem (\ref{prob1})  has infinite solutions.


\begin{thebibliography}{99}


\bibitem{adams}{ A. Adams and J. F. Fournier,}{\it \,  Sobolev spaces,  2nd ed.}, Academic Press, (2003).

\bibitem{adams2}{A. Adams and L.I. Hedberg,}{\it \, Function spaces and Potential Theory}{\it \, Grundlehren der Mathematischen Wissenschaften, vol. 314}, Springer-Verlag, Berlin, 1996.

\bibitem{Barreiro} {C.O. Alves and J.L.P. Barreiro,} {\it \,
Existence and multiplicity of solutions for a $p(x)$-Laplacian
equation with critical growth}, J. Math. Anal. Appl. (Print), 2013.

\bibitem{Abrantes} {C.O. Alves, J.V. Gon\c calves and J.A. Santos,} {\it Strongly Nonlinear Multivalued Elliptic
Equations on a Bounded Domain}, J Glob Optim, DOI
10.1007/s10898-013-0052-3, 2013.

\bibitem{BBR}{G. Bonano, G.M. Bisci and V. Radulescu,}{\it \, Quasilinear elliptic non-homogeneous Dirichlet problems through Orlicz-Sobolev spaces,} Nonl. Anal. 75 (2012), 4441-4456.

\bibitem{Brezis_Lieb} H. Brezis. and E. Lieb, {\it A relation between pointwise convergence of functions and
convergence of functinals}, Proc. Amer. Math. Soc. 88 (1983),
486-490.

\bibitem{Cerny} {R. Cern\'y,}{\it \, Generalized Moser-Trudinger inequality for unbounded domains and its application, } Nonlinear Differ. Equ. Appl. DOI 10.1007/s00030-011-0143-0.

\bibitem{chabro}
J. Chabrowski, Weak covergence methods for semilinear elliptic
  equations. World Scientific Publishing Company, 1999.

\bibitem{VGMS}{Ph. Cl\'ement, M. Garcia-Huidobro, R. Man\'asevich and K. Schmitt,}{\it \, Mountain pass type solutions for quasilinear elliptic equations,} Calc. Var. 11 (2000), 33-62.

\bibitem{Donaldson}{ T.K. Donaldson,}{\it \, Nonlinear elliptic boundary value problems in Orlicz-Sobolev spaces,} J. Diff. Equat. 10 (1971), 507-528.

\bibitem{Donaldson2}{T.K. Donaldson and N.S. Trudinger,}{\it \, Orlicz-Sobolev spaces and imbedding theorems,} J. Funct. Anal. 8 (1971) 52-75.


\bibitem{Fuchs1}{M. Fuchs and G. Li,} {\it Variational inequalities for energy functionals with nonstandard growth conditions}, Abstr. Appl. Anal. 3 (1998), 405-412.

\bibitem{Fuchs2}{M. Fuchs and V. Osmolovski,} {\it Variational integrals on Orlicz-Sobolev spaces. Z.} Anal. Anwendungen 17, 393-415 (1998) 6.

\bibitem{Fukagai1}{ N. Fukagai, M. Ito and K. Narukawa,}{\it \, Positive solutions of quasilinear elliptic equations with critical Orlicz-Sobolev nonlinearity on $\mathbb{R}^{N}$}, Funkcialaj Ekvacioj 49 (2006), 235-267.

\bibitem{Fukagai2}{ N. Fukagai, M. Ito and K. Narukawa,}{\it \, Quasilinear elliptic equations with slowly growing principal part and critical Orlicz-Sobolev nonlinear term}, Proc. R. S. Edinburgh 139 A (2009), 73-106.

\bibitem{Fukagai3}{ N. Fukagai and K. Narukawa,}{\it On the existence of multiple positive solutions of quasilinear elliptic eigenvalue
problems}, Annali di Matematica, 186, 539-564, 2007.

\bibitem{Gossez}{J.P. Gossez,}{\it \, Nonlinear elliptic boundary value problem for equations with rapidly (or slowly) increasing coefficients,} Trans. Am. Math. Soc. 190 (1974), 163-205.

\bibitem{Maso} Dal Maso, G. \& Murat, F., {\it Almost everywhere convergence of gradients of solutions to nonlinear elliptic
    systems}, Nonlinear Anal. 31 (1998), 405-412.

\bibitem{MR1}{M. Mihailescu and V. Radulescu,}{\it \, Nonhomogeneous Neumann problems in Orlicz-Sobolev spaces}, C.R. Acad. Sci. Paris, Ser. I 346 (2008), 401-406.


\bibitem{MR2}{M. Mihailescu and V. Radulescu,}{\it \, Existence and multiplicity of solutions for a quasilinear nonhomogeneous problems: An Orlicz-Sobolev space setting}, J. Math. Anal. Appl. 330 (2007), 416-432.

\bibitem{MD}{M. Mihailescu and D. Repovs,}{\it \, Multiple solutions for a nonlinear and non-homogeneous problems in Orlicz-Sobolev spaces,} Appl. Math. Comput. 217 (2011), 6624-6632.


 \bibitem{Pohozaev} Pohozaev, S. L., {\it Eingenfunctions for the
    equation $\Delta u+\lambda f(u)=0,$} Soviet Math. Dokl. 6 (1965)
    1408-1411.
\bibitem{Oneill}{R. O'Neill,}{\it \, Fractional integration in Orlicz spaces,} Trans. Amer. Math. Soc. 115 (1965) 300-328.

\bibitem{Rabinowitz1}
P.H. Rabinowitz, Minimax methods in critical point theory with
  applications to diferential equations. CBMS Reg. Conf. series in math 65, 1984.


\bibitem{xinmin}
W. Zhihui \& W. Xinmin, A multiplicity result for quasilinear elliptic
  equations involving critical sobolev exponents. Nonlinear Analysis, Theory,
  Methods and Applications 18 (1992), 559--567.
\end{thebibliography}
\end{document}